\documentstyle[12pt]{article}

\begin{document}
\def\s{\sigma}
\def\B{{\bf B}}
\def\C{{\bf C}}
\def\CC{{\bf C}^2}
\def\Cn{{\bf C}^n}
\def\L{\mbox{Lip }}
\def\b{\partial}
\begin{center}
{\large\bf PEAK POINTS FOR PSEUDOCONVEX DOMAINS: A SURVEY}

\mbox{}

ALAN NOELL

\end{center}

\vspace{.5in}

{\bf 1. Introduction}

Let $D$ be a  domain in $\Cn$ with smooth (that is, $C^\infty$)
boundary. If $0\leq \alpha\leq \infty$ we denote by $A^\alpha (D)$
the space of functions holomorphic on $D$ and of class $C^\alpha$
on $\overline{D}$. We write $A(D)$ for $A^0 (D)$ and $A^\omega
(D)$ for the space of functions holomorphic on a neighborhood of
$\overline{D}$. We say that a point $p\in\partial D$ is a peak
point relative to $D$ for $A^\alpha (D)$ if there exists a
function $f\in A^\alpha (D)$ so that $f(p)=1$ and $|f|<1$ on
$\overline{D}\setminus \{p\}$. We call $f$ a peak function. This
condition is clearly equivalent to the existence of a strong
support function, that is, a function $g\in A^\alpha (D)$ so that
$g(p)=0$ and $\mbox{Re }g>0$ on $\overline{D}\setminus \{p\}$. We
say that $p\in \b D$ is a local peak point for $A^\alpha(D)$ if
$p$ is a peak point for $A^\alpha(D\cap V)$ for some neighborhood
$V$ of $p$.

We want to determine whether boundary points are peak points.
Finding a peak function can be thought of as giving a quantitative
converse to the maximum modulus principle. Now observe that if $f$
is a peak function at $p$ relative to $D$ then $1/(1-f)$ is a
holomorphic function on $D$ with no holomorphic extension past
$p$. Thus if every boundary point of $D$ is a peak point then $D$
is a domain of holomorphy. In light of this fact, if we want to
determine whether boundary points are peak points it makes sense
to restrict our attention to domains of holomorphy. By the
solution of the Levi problem these are the (Levi) pseudoconvex
domains. Here is another observation: If there is a complex disc
in the boundary of $D$ then (by the maximum modulus principle) no
point in the relative interior of that disc can be a peak point
for $A(D)$.  A natural condition in this context is that $D$ be of
finite type at the point $p$ in question, in the sense of
D'Angelo: There is a finite upper bound on the order of contact of
nontrivial complex analytic varieties with the boundary at $p$.
The infimum of all such upper bounds is called the type at $p$. We
formulate  the major open question in terms of this notion.

\proclaim Question. If $D$ is a smooth bounded pseudoconvex domain
in $\Cn$ of finite type at $p\in\partial D$, is $p$ a peak point
for $A^\alpha(D)$ for some $\alpha>0$?

In this paper we survey the current knowledge with regard to this
question. The question has been answered in the affirmative in
$\CC$, but even here there are remaining problems, particularly in
connection with the relationship between the regularity exponent
$\alpha$ and the type at the boundary point in question. In our
exposition we present unpublished results of Laszlo that apply to
certain domains in $\CC$. His results give very good regularity
and also construct the peak function in a more explicit way than
was previously done. The question stated above is still open in
$\Cn$ for $n>2$. In fact, it is not known whether every boundary
point of a smooth bounded pseudoconvex domain $D$ of finite type
in $\Cn$ is a peak point for $A(D)$.

In Section 2 we provide some background information and analyze
the case of strict pseudoconvexity. Section 3 discusses the notion
of finite type, presents the Kohn-Nirenberg domain for which no
smooth peak function exists at a point, and considers notions of
strict type. In Section 4 we first outline the solution in $\CC$
by Bedford and Forn{\ae}ss of the main question above. Then we
present the modification  by Laszlo of the Bedford-Forn{\ae}ss
construction. Section 5 discusses other approaches to the problem,
and  Section 6 presents what is known about the question in $\Cn$.
The concluding section is a miscellany of results.

This work was supported in part by a Big 12 Faculty Fellowship
funded by Oklahoma State University. Thanks go to my host, the
Department of Mathematics at Texas A\&M University, and in
particular to the group in complex analysis there. Thanks also go
to Geza Laszlo for granting permission to present results from his
dissertation.

\vspace{.2in}

{\bf 2. Background information and the strictly pseudoconvex case}

\vspace{.2in}

Let $D$ be a domain in $\Cn$ with smooth boundary, and let $r$ be
a smooth defining function for $D$, in the following sense: For
some neighborhood $V$ of $\partial D$, the real-valued function
$r$ is smooth  on $V$, $$V\cap D=\{z\in V\colon r(z)<0\},$$ and
$\nabla r\neq 0$ on $V$. For $p\in\partial D$ and $t=(t_1,\ldots
,t_n)\in \Cn$ we write $\b r_p(t)$ for $$\sum_{j=1}^n \frac{\b
r}{\b z_j}(p) t_j \; , $$ and we write $L_p(r,t)$ for
$$\sum_{j,k=1}^n \frac{\b ^2 r}{\b z_j \b \bar{z}_k}(p) t_j
\bar{t}_k \; ,$$ the Levi form of $r$ at $p$ applied to $t$. The
domain $D$ is said to be (Levi) pseudoconvex if for each $p\in\b
D$ we have $L_p(r,t)\geq 0$ whenever $\b r_p(t)=0.$ The geometric
interpretation of the condition $\b r_p(t)=0$ is that $t$ belongs
to the complex tangent space to the boundary at $p$, that is, the
maximal complex subspace of the real tangent space at $p$.  The
domain is said to be strictly pseudoconvex at a boundary point $p$
if $L_p(r,t)> 0$ whenever $\b r_p(t)=0$ and $t\neq 0$. Note that
the complex tangent space to the boundary is trivial for domains
in ${\bf C}$ with smooth boundary, so such domains are strictly
pseudoconvex at each boundary point. In general the condition of
strict pseudoconvexity is independent of the choice of defining
function, and it implies that there exist a defining function
$\rho$ and a constant $c$ so that $$L_p(\rho,t)\geq c||t||^2 \; ,
t\in \Cn. \leqno(\ast) $$ In fact, if $r$ is a defining function
the choice $\rho=r+Mr^2$ satisfies condition $(\ast)$ for $M$
sufficiently large and $c$ sufficiently small.

This observation makes it easy to construct a local peak function
at a strictly pseudoconvex boundary point. The following result in
some form is usually attributed Levi, who introduced the notion of
pseudoconvexity we have presented (see \cite{lev1} and
\cite{lev2}).

\proclaim 2.1 Proposition. If $D$ is a smooth bounded domain in
$\Cn$ that is strictly pseudoconvex at $p\in\b D$ then $p$ is a
local peak point for $A^\omega(D)$.

\noindent{\bf Proof:} For ease of notation we assume that 0 is a
strictly pseudoconvex boundary point of $D$. Choose a defining
function $\rho$ for which condition $(\ast)$ holds. For
$w=(w_1,\ldots, w_n)\in \Cn$  put $$g(w)=-2\b
\rho_0(w)-\sum_{j,k=1}^n \frac{\b ^2 \rho}{\b z_j \b z_k}(0) w_j
w_k \; ,$$ a holomorphic polynomial vanishing at 0. Now the
complex form of the Taylor expansion of $\rho$ about 0 gives
$$\rho(z)=-\mbox{ Re }g(z)+L_0(\rho,z)+o(||z||^2).$$ But if
$z\in\overline{D}$ then $\rho(z)\leq 0$ and we have $$\mbox{ Re
}g(z)\geq L_0(\rho,z)+o(||z||^2).$$ Hence, by condition $(\ast)$,
$g$ is a local strong support function at 0. \hfill $\Box$

\vspace{.1in}

Note that the local strong support function is expressed in an
elementary way in terms of the second-order Taylor expansion of
the special defining function.

To get a global peak function we need the following result due to
Kohn (see \cite{ko:gl} and pages 229--231 of \cite{ko:gl2}).

\proclaim 2.2 Theorem. Let $D$ be a smooth bounded pseudoconvex
domain in $\Cn$ and $\phi$ a smooth $(0,1)$-form on
$\overline{D}$, in the sense that $\phi=\phi_1d\bar{z}_1+\cdots
+\phi_n d\bar{z}_n$ with $\phi_j\in C^\infty(\overline{D})$ for
all $j$. Assume that $\phi$ is $\overline{\b}$-closed, that is,
$\displaystyle\frac{\b \phi_j}{\b \bar{z}_k}=\frac{\b \phi_k}{\b
\bar{z}_j}$  for all $j,k$. Then there exists $\psi\in
C^\infty(\overline{D})$ so that $\overline{\b}\psi=\phi$, in the
sense that $\displaystyle\frac{\b \psi}{\b \bar{z}_j}=\phi_j$  for
all $j$.

Now we can show that  every strictly pseudoconvex boundary point
of a smooth bounded pseudoconvex domain is a (global) peak point.
This was proved by Hakim and Sibony (see \cite{hasib:fs},
\cite{hasib:qc}) and by Pflug in \cite{pf}.

\proclaim 2.3 Theorem. If $D$ is a smooth bounded pseudoconvex
domain in $\Cn$ that is strictly pseudoconvex at $p\in\b D$ then
$p$ is a peak point for $A^\infty(D)$.

\noindent{\bf Proof:} Apply Proposition 2.1 to get a local smooth
strong support function $g$ at $p$. Choose a nonnegative
$C^\infty$ function $\chi$ with support in a small neighborhood of
$p$ so that $\chi=1$ near $p$, and define the $(0,1)$-form $\phi$
on $\overline{D}$ by $\phi(p)=0$ and $\phi=\overline{\b} (\chi
/g)$ away from $p$. Then $\phi$ is smooth and
$\overline{\b}$-closed, so by Theorem 2.2 there exists $\psi\in
C^\infty(\overline{D})$ so that $\overline{\b}\psi=\phi.$ Adjust
$\psi$ by adding a constant to ensure that it has negative real
part. Now define $G$ on $\overline{D}$ by $G(p)=0$ and $G=1/(\chi
/g -\psi)$  away from $p$. It is easy to see that $G$ belongs to
$A^\infty(D)$ and is a strong support function at $p$ relative to
$D$. \hfill $\Box$

\proclaim 2.4 Remark.

\noindent 1) The method of argument in the proof of Theorem 2.3
can be adapted as follows. Let $p$ be a boundary point of a smooth
bounded pseudoconvex domain $D$. Suppose we are given a local peak
function $f$ at $p$ belonging to $A^\alpha(D)$ for some $\alpha
\leq \infty$. Assume that $f$ is smooth up to the boundary away
from $p$. Then $p$ is a (global) peak point for $A^\alpha(D)$.
Because of this fact, in the subsequent exposition we will focus
on constructing such local peak functions. Typically the local
peak function will extend to be holomorphic across the boundary
away from $p$. Sometimes the domain under consideration will  be
unbounded; Theorem 2.2 has been extended by Gay and Sebbar to
unbounded pseudoconvex domains with smooth boundary. See
Th\'eor\`eme 4.1 of \cite{gase}.

\mbox{}

\noindent 2) In Theorem 4.4 of \cite{ross}, Rossi proved that if
$\overline{D}$ has a neighborhood basis of pseudoconvex domains
then every local peak point for $A^\omega (D)$ is a peak point for
$A^\omega (D)$.
 If $D$ strictly
pseudoconvex at every boundary point then it is easy to see
        that $\overline{D}$ has a neighborhood basis of pseudoconvex
domains, so every boundary point is a peak point for $A^\omega
(D)$. This condition on the neighborhood basis also holds for
pseudoconvex domains of finite type, but the proof is much more
difficult. Rossi's result is not valid without this condition on
the neighborhood basis: In \cite{difo}, Diederich and Forn{\ae}ss
defined a family of smooth bounded pseudoconvex domains with no
such neighborhood basis. (These ``worm'' domains have a complex
annulus in the boundary.) In Theorem 3 of that paper they showed
that there are strictly pseudoconvex boundary points that are not
(global) peak points for $A^\omega $.

\mbox{}

\noindent 3) In \cite{sib:dbar}, Sibony constructed a smooth
bounded pseudoconvex domain $D$ in ${\bf C}^3$ having the
following property: There is a point $p\in\b D$ that is a local
peak point for $A(D)$ but is not a peak point. Relying on the
method of proof we described in Theorem 2.3, Sibony draws the
following conclusion: There exists a $(0,1)$-form $\phi$ on $D$
whose coefficients are smooth on $D$ and continuous on
$\overline{D}$ so that $\phi$ is $\overline{\b}$-closed but no
solution of $\overline{\b}\psi=\phi$ is bounded on $D$. (Compare
this with Theorem 2.2.)

\mbox{} \noindent

4) Let $D$ be a smooth bounded pseudoconvex domain. Assume that
$D$ is of finite type (as defined in Section 1), or more generally
that $\b D$ contains no nontrivial complex manifold. It turns out
that the set of boundary points at which $D$ is strictly
pseudoconvex is dense in the boundary of $D$. (If not, by
Freeman's work in \cite{free}
    there would be a local foliation of the boundary by complex
submanifolds, a contradiction. For a precise statement in the case
of finite type see the Main Theorem of Catlin's paper
\cite{catl}.) By Theorem 2.3, then, for pseudoconvex domains of
finite type the set of peak points for $A^\infty(D)$ is dense in
the boundary of $D$. Another result along these lines is due to
Basener: In \cite{ba} he showed that if $p$ is a $C^2$ boundary
point of an open set $D$ in $\Cn$ and $p$ is a peak point for
$A(D)$ then $p$ is a limit of strictly pseudoconvex boundary
points of $D$.

\vspace{.1in}

The following standard result is basic for studying peak functions
with some regularity. For a proof see, e.g., Proposition 12.2 of
\cite{fost}.

\proclaim 2.5 Hopf lemma. Let $\Omega$ be a bounded domain in the
complex plane with $C^2$ boundary and $u$ a negative subharmonic
function on $\Omega$. For $z\in \Omega$ denote by $\delta (z)$ the
distance from $z$ to $\b \Omega$. Then there exists a constant
$c>0$ so that $u(z)\leq -c\delta (z)$ for all $z\in \Omega$.

Now we apply the Hopf lemma to peak functions.

\proclaim 2.6 Holomorphic Hopf lemma. Let $D$ be an open set in
$\Cn$ with $C^2$ boundary near $p\in \b D$. Assume that $p$ is a
local peak point for $A^1(D)$ with strong support function $g$.
Then the derivative at $p$ of $\mbox{Re }g$ in the direction of
the outward normal is negative.

\noindent{\bf Proof:} Apply the Hopf lemma (2.5) to $-\mbox{Re }g$
near $p$  on the intersection of $D$ with the complex line spanned
by a normal to $\b D$ at $p$. \hfill $\Box$

\vspace{.1in}

We conclude this section with one version of the Bishop
$1/4$--$3/4$ method for constructing peak functions.

\proclaim 2.7 Theorem. Let $D$ be a bounded domain in $\Cn$ and
$p\in\b D$. Assume that for every neighborhood $V$ of $p$ there
exists $f\in A(D)$ so that $|f|\leq 1$ on $\overline{D}$, $f(p)>
3/4$, and $|f|< 1/4$ on $\overline{D}\setminus V$. Then $p$ is a
peak point for $A(D)$.

\vspace{.1in}

The method in fact applies much more generally to construct a peak
function from a family of approximate peak functions. The idea of
the proof is to construct inductively a nested sequence of
neighborhoods of $p$ and to add (scaled) approximate peak
functions for those neighborhoods. For a statement and proof in
the case of certain closed subspaces of $C(X)$ (with $X$ compact
Hausdorff), see Theorem 2.3.2 of \cite{brow}.

\vspace{.2in}

{\bf 3. Finite type, the Kohn-Nirenberg domain, and strict type}

\vspace{.2in}

Let $D$ be a smooth bounded pseudoconvex domain in $\Cn$. Recall
that $D$ is of finite type at a point $p\in\b D$ if there is a
finite upper bound on the order of contact of complex analytic
varieties with the boundary at $p$. The infimum of all such upper
bounds is the type at $p$. We remark that every bounded
pseudoconvex domain with real-analytic boundary is of finite type.
See D'Angelo's book \cite{dabo} for further background. In this
section we will re-formulate this notion for domains in $\CC$.  In
addition, we consider certain model domains having no smooth peak
function at a point. We also consider so-called strict type
conditions in $\Cn$.

First we make the following observation. Fix a domain $D$ in $\Cn$
with smooth boundary and a point $p\in \b D$.  Suppose we are
given local holomorphic coordinates $(z,w)$, with $z\in {\bf
C}^{n-1}$ and $w=u+iv$, in which $p=0$ and $u$ points in the
outward normal direction to $\b D$ at $p$. Then $D$ has a smooth
defining function of the form $$u+A(z)+B(z)v+O(|v|^2),$$ where $A$
vanishes to order at least 2 at 0 and $B(0)=0$. Note that $A$ and
$B$ depend on the choice of coordinates. In the rest of this paper
we refer to this expression as a standard form for a defining
function in these coordinates.

Now we restrict attention to pseudoconvex domains in $\CC$ and use
the preceding setup to  re-formulate the notion of finite type. If
in standard form for some coordinates we have $A\equiv 0$, then
the complex manifold $w=0$ lies in the boundary, so $D$ is not of
finite type at $p$. Now assume that $A\not\equiv 0$ and write
$$A(z)=P(z)+O(|z|^{m+1})$$ with $P\not\equiv 0$ a polynomial (in
$z\in {\bf C}$ and $\bar{z}$) that is homogeneous of degree $m\geq
2$. Note that $P$ depends on the choice of coordinates.  Then $D$
is of finite type at $p$ if and only if there exist such
coordinates in which $P$ is not harmonic. Further, the degree $m$
of such a $P$ is uniquely determined (independent of the choice of
coordinates and defining function) and is the type at $p$. For an
explanation of this re-formulation, see Lecture 28 of \cite{fost}.
It is easy to see that the pseudoconvexity of $D$ implies that $P$
is subharmonic. In particular, the type is an even integer.

The case of type 2 corresponds to strict pseudoconvexity, and each
such point is a local peak point for $A^\omega(D)$ and a global
peak point for $A^\infty(D)$ by the results of the preceding
section. It turns out that in $\CC$ points of type 4 have the same
property. We sketch the proof to indicate how similar results are
proved, but we leave details to be filled in by the indicated
references.

\proclaim 3.1 Proposition. Let $D$ be a smooth bounded
pseudoconvex domain in $\CC$. If $D$ is of type 4 at $p\in\b D$
then $p$ is a local peak point for $A^\omega(D)$.

\noindent{\bf Proof:} We use the preceding notation. Choose
coordinates so that $P$ is a subharmonic polynomial homogeneous of
degree 4. After a holomorphic change of coordinates we may assume
that $P$ has the form $$P(z)=|z|^4+t|z|^2\mbox{Re }(z^2)$$ with
$t\geq 0$. Because $P$ is subharmonic we have $t\leq 4/3$. Now the
key point is that if $\epsilon>0$ is sufficiently small then there
exists $c>0$ so that $$|z|^4+t|z|^2\mbox{Re
}(z^2)+\frac{t}{4(1-\epsilon)}\mbox{Re }(z^4)\geq c|z|^4.$$ (See
Lecture 28 of \cite{fost}.) Now replace $w$ by
$w+\displaystyle\frac{t}{4(1-\epsilon)}z^4$ for small
$\epsilon>0$. In the new coordinates the corresponding homogeneous
polynomial $P$ is the expression on the left-hand side of the
preceding inequality.  Now we analyze the function $B$. If the
leading term in the expansion of $B$ about 0 is linear, say
$\mbox{Im }(\beta z)$ for some $\beta\in {\bf C}$, replacing $w$
by $w(1+\beta z)$ results in a new expansion in which $B$ vanishes
to order at least 2 at 0. Further, the polynomial $P$ is not
changed. (See Lemma 2.8 and  the subsequent discussion in
\cite{bloo} or the proof of Proposition 1.1 in \cite{fosib}.)
Hence there exists a constant $T>0$ so that $$B^2(z)< TA(z)$$ if
$z\neq 0$ is sufficiently small. It follows that $-w-Mw^2$ is a
local strong support function if $M$ is a large constant. (See
Lemma 2.2 of \cite{bloo}.) \hfill $\Box$

\vspace{.1in}

The last part of this proof is a special case of a general result
due to Bloom. To explain this we use the standard form involving
the functions $A$ and $B$ developed for general domains.

\proclaim 3.2 Proposition. Let $D$ be a domain in $\Cn$ with
smooth boundary and $p\in \b D$. Then $p$ is a local peak point
for $A^\omega (D)$ if and only if the following holds: There exist
local holomorphic coordinates in which, for some $T>0$, $B^2(z)<
TA(z)$ if $z\neq 0$ is sufficiently small.

\vspace{.1in}

For the proof we refer the reader to Lemma 2.2 and Lemma 2.4 of
\cite{bloo}. We remark that in Bloom's characterization the local
peak function is given explicitly in terms of the special defining
function in local coordinates: The local strong support function
he constructs has the form $-w-Mw^2$ for some constant $M$,  as in
the proof of our Proposition 3.1.

\vspace{.1in}

We have seen that, for smooth bounded pseudoconvex domains in
$\CC$, points of type 2 and points of type 4 are local peak points
for $A^\omega$. For points of type 6 there is no such result. The
following is due to Forn{\ae}ss \cite{fo:c1}.

\proclaim 3.3 Proposition. Fix a number $t> 0$. Consider the
domain $D$ defined near the origin in $\CC$ by $$\mbox{Re
}w+|z|^6+t|z|^2\mbox{Re }(z^4)+|zw|^2<0.$$ If  $t\leq 9/5$ then
$D$ is pseudoconvex (and strictly pseudoconvex away from the
origin), and if $t>1$ the origin is not a local peak point for
$A^1(D)$.

\noindent{\bf Proof of weaker result:} We assume that $1<t\leq
9/5$. It is straightforward to verify the pseudoconvexity claims,
and we show only that there is no local strong support function at
the origin in $A^\omega$. (See the comments following the proof.)
Assume for a contradiction that such a function $g$ exists. Note
that if $(z,w)\in \overline{D}$ then $(iz,w)\in \overline{D}$.
Hence by averaging we may assume that $g(iz,w)\equiv g(z,w)$. Now
by the holomorphic Hopf lemma (2.6) and the implicit function
theorem, near the origin we can write the zero set of $g$ as a
graph $w=h(z)$ for some function $h$ holomorphic near 0 satisfying
$h(0)=0$. For ease of notation put $$P(z)= |z|^6+t|z|^2\mbox{Re
}(z^4).$$ Note that $P$ takes negative values in every
neighborhood of 0 because $t>1$.
 Now the zero set of $g$ does not intersect $D$ near the origin, and using the
defining function for $D$ gives $$\mbox{Re
}h(z)+P(z)+|zh(z)|^2\geq 0$$ for $z$ near 0. This shows that
$h\not\equiv 0$ because $P$ takes negative values. By the symmetry
of $g$ we have $h(iz)\equiv h(z)$, so the leading term of the
Taylor expansion of $h$ about 0 is $cz^{4k}$ for some positive
integer $k$ and constant $c$. But the preceding inequality and the
fact that $P$ takes negative values shows that $k< 2$, and because
$\mbox{Re }(cz^{4})$ takes negative values we have a
contradiction. \hfill $\Box$

\vspace{.1in}

A careful examination of the preceding proof will show that the
origin is not a local peak point for  $A^6$. (See also page 402 of
\cite{hasib:qc}.) To prove that there is no peak function in
$A^1$, Forn{\ae}ss in \cite{fo:c1} starts with the Hopf lemma as
above and obtains a normal family of harmonic functions by
examining the (normalized) values of a strong support function on
slices of the domain by the planes $w=-\epsilon$.

Shortly after Forn{\ae}ss proved this result, in \cite{bf1}
Bedford and he constructed, at each boundary point of a
pseudoconvex domain of finite type  in $\CC$, a peak function in
$A^\alpha$ for some $\alpha >0$. In the next section we will
discuss this result and also see from Laszlo's work that for the
domain in Proposition 3.3 (in the case $t=9/5$) we may take
$\alpha=1/18$. It is an open problem to determine for which values
of $\alpha <1$ there is a peak function at the origin in
$A^\alpha$ relative to this domain. The current methods seem
incapable of producing a peak function in $A^\alpha$ for $\alpha$
greater than $1/6$, the reciprocal of the type. In particular, it
is not known whether the origin is a peak point for $A^\alpha$
with $\alpha <1$ arbitrarily close to 1.

It is clear that certain properties of the homogeneous polynomial
$$P(z)= |z|^6+t|z|^2\mbox{Re }(z^4)$$ appearing in the defining
function are crucial for the proof of Proposition 3.3. This
polynomial is a special case (with $k=3$ and $\ell=2$) of a more
general two-term homogeneous subharmonic polynomial, namely
$$P(z)=|z|^{2k}+t|z|^{2k-2\ell}\mbox{Re }(z^{2\ell}),$$ where $k$
and  $\ell$ are positive integers so that $\ell<k$, and $0<t\leq
k^2/(k^2-\ell^2)$. Using such a polynomial to define a domain as
in Proposition 3.3 gives a pseudoconvex domain of type $2k$ at the
origin. The seminal study of such a domain was the paper
\cite{koni} by Kohn and Nirenberg, which considered the case $k=4$
and $\ell=3$ with $t=15/7$. They showed that if the zero set of a
function holomorphic near the origin passes through the origin
then the zero set intersects the domain. (See also Proposition 1.2
of  \cite{sib:asp}.) This was quite an unexpected result;
apparently the expectation in the early 1970s was that some sort
of convexity condition would carry over from the strictly
pseudoconvex case.

We remark that for such a non-convexity result one needs not only
that $t>1$ but also that $\ell$ does not divide $k$ (consider the
case of type 4). In fact, if $\ell$ divides $k$ then the origin is
a local peak point for $A^\omega$. Kol\'a\v{r} proves this and
several related results in \cite{kol}.

In light of the Kohn-Nirenberg domain there are two possible
approaches to studying the existence of peak functions. One is to
try to construct peak functions with minimal regularity on domains
of finite type. This was the approach of Bedford and Forn{\ae}ss,
as mentioned above, and this will be the point of view in much of
the remainder of our study. Another approach is to add conditions
stronger than finite type that guarantee the existence of more
regular peak functions. These are often called strict type
conditions, and in the remainder of this section we consider them.

In \cite{ko:st}, Kohn introduced a strict type condition for
domains in $\CC$ and showed that for smooth bounded pseudoconvex
domains in $\CC$ every point of strict type (in that sense) is a
local peak point for $A^\omega$. In \cite{range1}, \cite{range2},
and \cite{range3}, Range gave a condition on domains in $\Cn$ that
in $\CC$ is weaker  than Kohn's condition. In terms of the
function $A$ appearing in the standard form for a defining
function, Range's condition is that there exist holomorphic
coordinates so that, for some constant $c>0$ and some integer $N$,
$$A(z)\geq c|z|^N$$ near $z=0$. In \cite{hasib:qc}, Hakim and
Sibony studied this strict type condition as it relates to peak
functions. In \cite{bloo}, Bloom generalized Kohn's condition to
$\Cn$ and introduced a third notion (weaker than Kohn's and
stronger than Range's) of strict type. That paper gives a very
helpful discussion of all three conditions. To avoid
technicalities we give only the main ideas.

First we state one main result of Hakim and Sibony from
\cite{hasib:qc}.  The proof is omitted.

\proclaim 3.4 Theorem. Let $D$ be a smooth bounded domain (not
necessarily pseudoconvex) in $\Cn$ and $p$ a boundary point of
strict type in the sense of Range: There exist holomorphic
coordinates near $p$ so that, for some constant $c>0$ and some
integer $N$, $$A(z)\geq c|z|^N$$ near $z=0$. Then $p$ is a local
peak point for $A^\alpha (D)$ for some $\alpha >1$. In addition,
if we may take $N\leq 4$ then $p$ is a local peak point for
$A^\omega$.

\vspace{.1in}

As we noted above, the strict type condition introduced by Bloom
in \cite{bloo} is stronger than Range's condition. It, too,
involves a positive-definite condition, but one expressed in terms
of the leading homogeneous part of $A$ in weighted holomorphic
coordinates. Bloom proves the following result.

\proclaim 3.5 Theorem. Let $D$ be a smooth bounded pseudoconvex
domain and $p$ a boundary point of strict type in the sense of
Bloom. Then $p$ is a local peak point for $A^\omega$.

\noindent{\bf Sketch of proof:} The proof uses the result we
stated as Proposition 3.2: The goal is to verify the inequality
that, for some $T>0$, we have $B^2(z)< TA(z)$ if $z\neq 0$ is
sufficiently small. Using successive holomorphic changes of
coordinates allows one to ensure that the leading part of $B$
(that is, the part of lowest weight) is not the real part of a
holomorphic function. The pseudoconvexity condition then implies
that this part is of weight at least half that of the leading part
of $A$. Using  the positive-definite condition gives the desired
inequality. \hfill $\Box$

\vspace{.1in}

In \cite{bloo}, Bloom also provided a striking example. Among
other things,  this example shows that the notion of strict type
given by Range does not in general imply the existence of a local
peak function in $A^\omega$.

\proclaim 3.6 Proposition. Consider the domain $D$ defined near
the origin in $\CC$ by $r<0$, where
$$r(z,w)=u+100(|z|^{10}+|z|^2\mbox{Re }(z^8))+|z|^6v+|z|^2v^2.$$
(Here we use $(z,w)$ for coordinates and write $w=u+iv$.) Then $D$
is pseudoconvex (and strictly pseudoconvex away from the origin),
and the origin is not a local peak point for $A^\omega(D)$---in
fact, is not a local peak point for $A^{13}(D)$. The origin is,
however, of strict type in the sense of Range, and so is a local
peak point for $A^\alpha (D)$ for some $\alpha
>1$.

\noindent{\bf Sketch of proof:}  It is straightforward to verify
the pseudoconvexity claims. Bloom shows that there are no local
coordinates in which the basic inequality $B^2(z)< TA(z)$ (for
$z\neq 0$ sufficiently small) holds. By Proposition 3.2 this
proves that there is no local peak function in $A^\omega(D)$. For
the rest, one checks that $r(z,z^{16})\geq c|z|^{16}$ for some
$c>0$ (see \cite{bloo} and Lecture 30 of \cite{fost}). This
implies that the origin is of strict type in the sense of Range,
and the concluding statement of the proposition follows from
Theorem 3.4. \hfill $\Box$

\vspace{.1in}

We remark that for this example it is not known for which $\alpha$
between 1 and 13 the origin is a local peak point for $A^\alpha$.

When Range introduced his notion of strict type, he observed that
this condition holds at all boundary points of bounded convex
domains with real-analytic boundary. We remark that clearly every
such boundary point is a peak point for $A^\omega$: Because the
boundary can contain no line segment, at each boundary point the
tangent plane intersects the closure of the domain only at the
point of tangency, so there is an affine strong support function.
A main motivation for  Range's work was to prove H\"older
estimates for the $\overline{\partial} $-equation, and for that
study he needed a smoothly varying family of functions satisfying
certain estimates. In more recent work, also motivated by the
study of the $\overline{\partial} $-equation,
 Diederich
and Forn{\ae}ss studied smooth bounded convex domains of finite
type. In \cite{difoc1}, they gave an explicit formula for a local
strong support function in $A^\omega$ at each boundary point of a
convex domain of finite type. In fact, they constructed a smoothly
varying family of such functions satisfying optimal estimates with
regard to the order of contact. See also their work in
\cite{difoc2}, where they extended these results to cover
so-called lineally convex domains of finite type. (These are
domains of finite type so that through each boundary point there
passes a complex hyperplane that does not intersect the domain.)
It follows that if a smooth bounded domain of finite type is
convex (or, more generally, lineally convex), then every boundary
point is a global peak point for $A^\omega$.

In part 4) of Remark 2.4 we noted that, for smooth bounded
pseudoconvex domains of finite type, the set of peak points for
$A^\infty(D)$ is dense in the boundary. Now we give a refinement
of that statement for bounded pseudoconvex domains with
real-analytic boundary  in $\CC$.  It follows from this refinement
that the points where the Kohn-Nirenberg behavior occurs form a
relatively small subset of the boundary.

The following theorem is drawn from the work of Noell and
Stens{\o}nes in \cite{no:st}. (See Proposition 1.1 and Proposition
1.6 there.)

\proclaim 3.7 Theorem. Let $D$ be a bounded pseudoconvex domain
with real-analytic boundary in $\CC$. Then there exists a set
$E\subset \b D $ that is a finite union of singleton sets and
real-analytic curves so that every point in $\b D\setminus E$ is a
peak point for $A^\infty (D)$.

\vspace{.1in}

Instead of giving the proof we explain the main point. We say that
a curve in the boundary is complex tangential at a point if its
tangent line lies in the complex tangent space to the boundary at
that point. The key fact is that the Kohn-Nirenberg phenomenon
cannot occur at any point belonging to a curve having both of the
following properties: the curve is complex-tangential at every
point, and the boundary is of constant finite type along the
curve. (The curves in the preceding theorem can be taken to be
nowhere complex-tangential.) For example, the following is a
consequence of Noell's work in \cite{no:ct}. (In \cite{no:ct} the
conclusion is that the curve is locally a peak set for $A^\infty
(D)$, but by examining the proof or using Theorem 4.10 of
\cite{no:peak} one sees that each point on the curve is a peak
point. The notion of peak set is defined in Section 7 below.)

\proclaim 3.8 Theorem. Let $D$ be a smooth bounded pseudoconvex
domain in $\CC$ and $\gamma$ a smooth curve in $\b D$ that is
complex tangential at every point. Assume that $D$ is of constant
finite type along $\gamma$. Then each point of $\gamma$ is a peak
point for $A^\infty (D)$.

To conclude this section we note Iordan's work in \cite{iord}
giving a sufficient condition for every boundary point to be a
peak point for $A^\infty$.

\vspace{.2in}

{\bf 4. The sector method}

\vspace{.2in}

As we noted earlier, Bedford and Forn{\ae}ss  in \cite{bf1}
obtained a positive answer to our main question for domains in
$\CC$, as follows.

\proclaim 4.1 Theorem. Let $D$ be a smooth bounded pseudoconvex
domain in $\CC$ and $p\in\partial D$ a point of finite type. Then
there exists $\alpha >0$ such that $p$ is a peak point for
$A^\alpha (D)$.

In this section we discuss their method for constructing peak
functions. They blow up the boundary point $p$ to obtain a line
bundle over a Riemann surface whose points are sectors in a cone
naturally associated to $D$ at $p$. Using an existence theorem
they produce (in the abstract) a section of this line bundle, and
the existence of a peak function follows. Here we also present
Laszlo's modification of the sector method. For certain domains he
obtains the desired section concretely, in terms of a holomorphic
polynomial naturally associated to the defining function. Because
the results in his dissertation \cite{la} have never been
published, we give several details of his method (but omit the
more technical parts).

We recall the discussion at the beginning of Section 3: If $p$ is
a point of finite type in the boundary of a smooth bounded
pseudoconvex domain $D$ in $\CC$, then there exist local
holomorphic coordinates $(z,w)$, with $w=u+iv$, in which $p=0$ and
$D$ has a smooth defining function of the form
$$u+P(z)+O(|v|^2+|zv|+|z|^{2k+1}).$$ Here $P$ is a subharmonic,
but not harmonic, polynomial on ${\bf C}$ that is homogeneous of
degree $2k$, the type of $D$ at $p$. For the time being we fix
such a polynomial $P$. The method of construction in \cite{bf1}
produces a peak function at $0$ relative to the domain
$$D'=\{(z,w)\in\CC\colon \mbox{Re }w+P(z)<0\}$$ by constructing an
associated function on the open cone $$C=\{(z,w)\in\CC\colon
\mbox{Re }(w^{2k})+P(z)<0\}.$$ Because of a bumping procedure the
associated function is in fact defined on a conical neighborhood
of $C$. (We postpone the description of such a bumping until
Section 6.) This means that the peak function relative to $D'$ is
actually defined on a neighborhood of $\overline{D}\setminus\{p\}$
near $p$: By bumping $C$ we in effect absorb the terms of higher
order in the special defining function for $D$. We thus obtain a
local peak function at $p$ relative to $D$. In what follows we
will drop all reference to the neighborhoods and work just on $D'$
and $C$.

Now we describe in more detail the construction in \cite{bf1} of a
peak function at $0$ relative to $D'$. As we stated above, the
peak function is defined in terms of an associated function that
is holomorphic on $C$. Here is a key observation: Intersecting $C$
with a complex line through the origin results in a finite
collection of sectors in that line. The associated function
constructed  in \cite{bf1} is linear on each of these sectors.
Given such a function $g$, one can define
$$h(z,w)=\prod_{j=1}^{2k} g(z,e^{i\pi j/k}w)$$ and then define a
function $H$ holomorphic on $D'$   by $$H(z,w^{2k})=h(z,w)$$ and
$H(0,0)=0$. If we take a large enough root of $H$ we obtain a
strong support function relative to $D'$. Because we are concerned
with the regularity at 0, we note that the function $g$, it turns
out, satisfies a H\"older condition of order 1, so $H$ is of class
$C^{1/(2k)}$. Then the strong support function belongs to
$A^\alpha$ for some $\alpha
>0$, but we have no estimate for $\alpha$ without information
about how $g$ winds around the origin.

The main problem, then, is to construct the function $g$ that is
linear on each sector. The difficulty is that the sectors may vary
with the complex line in a complicated way; for example, if we
follow a given initial sector around a loop in complex projective
space we might end up at a different sector. A key insight of
\cite{bf1} is that the sectors can be considered as points of a
Riemann surface over complex projective space. The function $g$ is
then obtained from a section of a line bundle over this Riemann
surface, as we noted in our introduction to the sector method.
Among other things, this construction requires showing that the
sectors vary in a smooth way with the complex line, a point we
address below. In the end the existence of the section follows by
solving a multiplicative Cousin problem.

This is an ingenious method, and the analysis of the domain in
terms of these sectors has proven useful in various contexts (for
example, \cite{bf2}, \cite{fo:dbar}, \cite{fore}, and
\cite{no:st}). One drawback is that the resulting peak function is
produced by means of an abstract existence theorem, which yields
little insight on how the peak function depends on the defining
function. The method also gives little control on the regularity
of the peak function. Laszlo's modification in \cite{la} of the
Bedford-Forn{\ae}ss construction addresses both of these concerns:
The peak function is defined in terms of the roots of an algebraic
equation naturally associated to the homogeneous polynomial $P$,
and there is good control on the regularity of the peak function.
His method applies when $P$ has a special form. Before we describe
this modified construction we take a closer look at how the
sectors vary with the complex line.

In a sequence of lemmas, Bedford and Forn{\ae}ss show how the
subharmonicity of $P$ is reflected in the behavior of the sectors.
The key point is that each sector has an angular opening of size
at most $\pi/(2k)$, and any two sectors must be separated by an
angular opening of size at least $\pi/(2k)$. These facts imply
that a sector cannot split as the complex line varies, a crucial
observation in showing that the sectors vary smoothly.

Here is the main result proved by Laszlo in \cite{la}.

\proclaim 4.2 Theorem. Fix positive integers $k$ and  $\ell$ so
that $k/2\leq \ell<k$, and put $t=k^2/(k^2-\ell^2)$. Define
$$P(z)=|z|^{2k}+t|z|^{2k-2\ell}\mbox{Re }(z^{2\ell}),$$ a
homogeneous subharmonic polynomial. Then 0 is a peak point
relative to every smooth bounded pseudoconvex domain in $\CC$
having a defining function of the form $$\mbox{Re
}w+P(z)+O(|w|^2+|zw|+|z|^{2k+1}).$$ Further, the peak function
constructed belongs to $A^\alpha$ with $\alpha=1/(6k)$. Note that
$6k$ is 3 times the type of the domain at 0.

We will see that the peak function is defined in terms of the
zeros of a holomorphic polynomial naturally associated to the
homogeneous polynomial $P$.

We sketch the proof of the theorem in two parts: the construction
of the Riemann surface, which is more concrete than the original
construction, and the definition of the function $g$ that is
linear on each sector. The first part is valid for  any polynomial
$P$ on ${\bf C}$ that is subharmonic, but not harmonic, and
homogeneous of degree $2k$. Hence we describe the construction for
any such polynomial. The second part applies only to the special
two-term polynomials defined in the theorem.

\vspace{.1in}

\noindent{\bf Construction of the Riemann surface:}  We write $P$
in the form $$P(z)=\sum_{j=0}^{2k} c_jz^j\bar{z}^{2k-j}.$$  Note
that $\bar{c_j}=c_{2k-j}$ for $0\leq j\leq 2k$ because $P$ is
real-valued, and $c_k>0$ because $P$ is subharmonic but not
harmonic. We assume that the harmonic terms in $P$ have been
removed, so $c_0=c_{2k}=0$. Put $a_j=c_{k+j}$ for $0\leq j\leq
k-1$, and for $z\in\C$ define $$Q_0(z)=\sum_{j=0}^{k-1}
a_jz^{2j}.$$ We think of this as the polynomial corresponding to
intersecting the cone $C$ with the complex line $w=0$. More
generally, for fixed $\zeta\in\C$ the intersection of $C$ with the
complex line $w=\zeta z$ is described by  $$\mbox{Re }(\zeta^{2k}
z^{2k})+P(z)<0,$$ and we define $$Q_\zeta(z)=\zeta^{2k}
z^{2k}+\sum_{j=0}^{k-1} a_jz^{2j}=\zeta^{2k} z^{2k}+Q_0(z).$$ This
is the holomorphic polynomial whose zeros will be used to define
the locally linear function $g$---the function from which the peak
function is derived.

The connection between $Q_\zeta(z)$ and the defining function can
be explained geometrically as follows. It is easy to see that, if
we write $z=|z|e^{i\theta}$ for $z\neq 0$, then we have the
equation $$\mbox{Re }(\zeta^{2k} z^{2k})+P(z)=|z|^{2k}\mbox{ Re
}Q_\zeta(e^{i\theta}). \leqno(\ast) $$
    Now the sectors in the complex line
$w=\zeta z$ divide the unit circle $|z|=1$ into a collection of
arcs, and equation ($\ast$) says that these arcs are the
components of the intersection with that circle of the set
$$E_\zeta=\{z\in\C\colon \mbox{ Re }Q_\zeta(z)<0\}.$$

Now instead of studying the sectors we study $E_\zeta$. The
intersection of $E_\zeta$ with the closed unit disc
$\overline{\B}$ has a relatively simple structure because $P$ is
subharmonic. We now describe this structure.

 A straightforward computation reveals
that the condition that $P$ be subharmonic is equivalent to
$$\mbox{Re }[4k^2Q_\zeta(z)-zQ_\zeta'(z)-z^2Q_\zeta''(z)]\geq 0
\mbox{ for } |z|=1. $$ Here we use the prime notation to indicate
the derivative with respect to $z$. Because the function on the
left is harmonic on the plane and nonzero at 0 (since
$a_0=c_k>0$), this condition in turn is equivalent to  $$\mbox{Re
}[4k^2Q_\zeta(z)-zQ_\zeta'(z)-z^2Q_\zeta''(z)]> 0 \mbox{ for }
|z|< 1. $$ It turns out that this condition implies the following:
Fix $\zeta$ and assume that $E_\zeta\cap\B$ is nonempty. Then
there exists $m$ with $1\leq m\leq 2k$ so that we have the
disjoint unions $$E_\zeta\cap\B=\bigcup_{j=1}^m W_j,$$
$$E_\zeta\cap\partial\B=\bigcup_{j=1}^m A_j,$$ and
$$\{z\in\C\colon \mbox{ Re
}Q_\zeta(z)=0\}\cap\overline{\B}=\bigcup_{j=1}^m B_j.$$ Here for
each $j$ we have that $W_j$ is open in $\C$, $A_j$ is an open arc,
$B_j$ is a smooth curve, and $\partial W_j=A_j\cup B_j.$ Further,
we have the following fact: For each $j$ there is a unique point
$\eta_j$ on $B_j$ closest to the origin, and $\eta_j/|\eta_j|\in
A_j$. Clearly, by our geometric interpretation of equation
($\ast$) above, there is a one-to-one correspondence between
$\{A_j\}_{j=1}^m$ and the sectors, and we denote by $S_j$ the
sector corresponding to $A_j$. Then $\eta_j\in S_j$.  We call
$\eta_j$ the distinguished point of $S_j$. Because any two sectors
must be separated by an angular opening of size at least
$\pi/(2k)$, the distinguished points are separated from each other
by an angle of at least $\pi/(2k)$. This property is crucial for
constructing the Riemann surface from the sectors. We remark that
we may arrange (by a change of variables) that when $\zeta=0$ the
number of sectors is $2k$.

Now it is not hard to see that the distinguished points $\eta_j$
are the solutions of the system of equations $$\mbox{ Re
}Q_\zeta(\eta)=\mbox{ Im }[\eta Q_\zeta'(\eta)]=0, \; \eta\in
\B.$$ If we think of this as a system of equations in $\zeta$ and
$\eta$, it defines a Riemann surface with projection
$\pi(\zeta,\eta)=\zeta$. Here we use the fact that the
distinguished points are separated, so each $\eta_j$ can be viewed
as a local real-analytic diffeomorphism.

Laszlo also showed how to extend the Riemann surface by
considering the limit as $\zeta$ tends to $\infty$. This concludes
the construction of the Riemann surface. \hfill $\Box$

\vspace{.1in}

\noindent{\bf Construction of the locally linear function $g$:}
Now we turn to the construction of a holomorphic function $g$ that
is linear on each sector in the cone $C$. This depends on locating
certain zeros of $Q_\zeta$, and the details were carried out by
Laszlo in case $P$ has the form described in the theorem. Carrying
out the details requires, among other things, careful analysis of
certain trigonometric polynomials. The special form of $P$ makes
this analysis possible (but still somewhat tedious).

As before, we work in the complex line $w=\zeta z$. The idea is
that a given sector $S_j$ has a distinguished point $\eta_j$, and
we associate to $\eta_j$ a zero $\alpha_j$ of $Q_\zeta$. Recall
that $\eta_j$ belongs to $$\{z\in\C\colon \mbox{ Re
}Q_\zeta(z)=0\},$$ so we are requiring that $\alpha_j$ belong to
the intersection of that level set with the ``conjugate'' level
set $$\{z\in\C\colon \mbox{ Im }Q_\zeta(z)=0\}.$$

The desired regularity of the peak function at 0 will hold if the
argument of $z$ and $\alpha_j$ are close enough when $z$ belongs
to $S_j$. (Shortly we will explain this point in more detail.)
Laszlo constructs a correspondence between the distinguished point
$\eta_j$ for each sector and a zero $\alpha_j$ in a way that gives
a real-analytic diffeomorphism as the complex line varies, with
the additional property that
$$|\arg{\alpha_j}-\arg{\eta_j}|<\frac{\pi}{4k}.$$

Now we conclude the proof of the theorem with a  focus on how the
regularity of the peak function is proved. Pick $(z,w)\in C$, and
for simplicity assume that $z\neq 0$. Write $w=\zeta z$. The
intersection with $C$ of the corresponding complex line  is a
disjoint union of sectors $S_1, \ldots, S_m$. Thus $z$ belongs to
$S_j$ for some $j$. If (in the preceding notation) $\eta_j$ is the
distinguished point of $S_j$ and $\alpha_j$ is the corresponding
zero of $Q_\zeta$, define $$g(z,w)=\frac{z}{\alpha_j}.$$

The function $g$ is holomorphic on $C$, is linear on each sector,
and satisfies a H\"older condition of order 1 near 0. Further,
because we always have
$$|\arg{\alpha_j}-\arg{\eta_j}|<\frac{\pi}{4k},$$ and because each
sector has an angular opening of size at most $\pi/(2k)$, we have
$$|\arg{g(z,w)}|<\frac{3\pi}{4k}.$$ Now define
$$h(z,w)=\prod_{j=1}^{2k} g(z,e^{i\pi j/k}w),$$  and define a
function $H$ holomorphic on $D'$ by $$H(z,w^{2k})=h(z,w)$$ and
$H(0,0)=0$. Then $H$ is of class $C^{ 1/(2k)}$ near 0. Now an
examination of the equation defining $C$ (in which $w$ appears to
the power $2k$) and the construction of $g$ shows that if
$w_0^{2k}=w$ then $$H(z,w)=[g(z,w_0)]^{2k}.$$ Hence
$$|\arg{H(z,w)}|<\frac{3\pi}{2}.$$ Now the argument of $H$ is
continuous, and we conclude that $H$ has a cube root that is a
strong support function and belongs to $A^{ 1/(6k)}$ near 0.
\hfill $\Box$

%\vspace{.1in}

\vspace{.2in}

{\bf 5. Alternative constructions of peak functions}

\vspace{.2in}

In this section we present two constructions of peak functions
that are alternatives to the sector method presented in the
preceding section. Both use an adaptation of the Bishop
$1/4$--$3/4$ method (see Theorem 2.7).

In \cite{fosib}, Forn{\ae}ss and Sibony constructed peak functions
on pseudoconvex domains of finite type in $\CC$. A key element of
their proof is an existence theorem for entire functions dominated
in terms of subharmonic functions, as follows.

\proclaim 5.1 Theorem. Let $\phi$ be a subharmonic function on
{\bf C} that is not harmonic. Assume that there exist constants
$m$ and $C$ so that $$\phi(z+z')-\phi(z)\leq C \mbox{ when }
|z'|\leq(1+|z|)^{-m}.$$ Then for each sufficiently large $\lambda$
and $M$ there exists an entire function $f$ so that $f(0)=1$ and
$$|f(z)|\leq M\frac{\exp(\lambda\phi(z))}{1+|z|} \mbox{ for } z\in
{\bf C}. $$

\vspace{.1in}

The proof depends on H{\"o}rmander's theory of solving the
$\overline{\partial} $-equation in $L^2$ spaces with weights. We
omit the details.

Now we can explain  the approach of \cite{fosib} to constructing
peak functions at points of finite type.

\proclaim 5.2 Theorem. Let $D$ be a smooth bounded pseudoconvex
domain in $\CC$ and $p\in\b D$ a point of finite type. Then $p$ is
a peak point for $A(D)$.

\noindent{\bf Sketch of proof:} A careful analysis of the defining
function in special coordinates shows that the domain can be
``bumped out'' away from $p$. (We will explain such bumping
methods in the next section.) The result is that it suffices to
construct a peak function at $p=0$ relative to the domain $U$
defined by $\mbox{Re }w+\phi(z)<0$. Here $\phi$ is a subharmonic
function that is homogeneous of degree $2k$ (the type of $D$ at
$p$) and of class $C^\infty$ away from 0. Now apply Theorem 5.1 to
$\phi$ with $m=2k$ to get an entire function $f$ for some choice
of $\lambda$ and $M$. If we put, for $j$ a positive integer,
$$f_j(z,w)=\exp(jw)f(j^{1/(2k)}\lambda^{-1/(2k)}z),$$ then
$f_j(0,0)=1$, and for $(z,w)\in \overline{U}$ we have
$$|f_j(z,w)|\leq M/(1+j^{1/(2k)}\lambda^{-1/(2k)}|z|).$$ Hence for
each $\delta>0$ and $\epsilon>0$ there exists $J$ so that if
$j\geq J$ and $|z|\geq\delta$  then $|f_j(z,w)|<\epsilon$ for all
$w$. An adaptation of the Bishop $1/4$--$3/4$ method  gives a peak
function for $p$ in $A(U)$. \hfill $\Box$

\vspace{.1in}

Another approach to constructing peak functions on domains of
finite type in $\CC$ is due to Forn{\ae}ss and McNeal in
\cite{fomc}. Their method is especially promising because it
produces local peak functions on domains satisfying certain
general conditions related to the Bergman kernel and the
$\overline{\partial} $-Neumann operator, both of which have been
objects of intense study.  A careful statement of these conditions
would require considerable background material, so we only sketch
the main ideas.

First we recall the definition of the Bergman kernel. Let $D$ be a
bounded domain in $\Cn$, and consider the space ${\cal H}^2(D)$ of
holomorphic functions belonging to $L^2(D)$, equipped with the
standard inner product. For fixed $q\in D$ the bounded linear
functional sending $f\in {\cal H}^2(D)$ to $f(q)$ can be
represented as the inner product of $f$ with a function $k_q\in
{\cal H}^2(D)$. The Bergman kernel for $D$ is the function $K$ on
$D\times D$ defined by $K(z,q)=k_q(z)$. It is a fact that always
$K(z,z)>0$ for $z\in D$.

Here is how the properties of the Bergman kernel are used in
\cite{fomc} to construct approximate peak functions.

\proclaim 5.3 Lemma. Let $D$ be a smooth bounded pseudoconvex
domain in $\CC$ of finite type and $p\in\b D$. Denote by $K$ the
Bergman kernel for $D$.  There exists a constant $C$ so that for
every neighborhood $V$ of $p$ the following holds: For $q$
sufficiently close to $p$ on the inward normal to $\b D$ at $p$,
define $f_q$ on $D$ by $$f_q(z)=K(z,q)/K(q,q).$$ Then $f_q$ is
holomorphic on $D$, $f_q(q)=1$,   $|f_q|\leq C$ on $D$, and
$|f_q|<1/2$ on $D\setminus V$.

\noindent{\bf Idea of proof:} Only the last two properties require
verification. For the last of these the key properties of the
Bergman kernel are the following. First, $K(z,z)$ tends to
infinity as $z$ approaches $\b D$; this holds (with good control
on the rate) for all smooth bounded pseudoconvex domains by the
work of Pflug in  \cite{pf2}. This property gives a good lower
bound for the denominator of $f_q$ in terms of the distance to the
boundary. Second, the Bergman kernel decays away from the boundary
diagonal; this depends (via work of McNeal in \cite{mc}) on the
existence of a subelliptic estimate for the $\overline{\partial}
$-Neumann problem, as proved by Catlin for pseudoconvex domains of
finite type in $\Cn$. Combining these two properties shows that we
can make $|f_q|$ small on $D\setminus V$ by taking $q$ close to
$p$. The remaining claim in the lemma is the uniform boundedness
on $D$ of the family $\{f_q\}$. This requires delicate estimates
of the kernel on (biholomorphic images of) polydiscs adapted to
the boundary geometry; such estimates are known to hold on domains
of finite type in $\CC$ and on certain other domains (see Section
6 below). \hfill $\Box$

\vspace{.1in}

With this lemma in hand, Forn{\ae}ss and McNeal adapt the Bishop
$1/4$--$3/4$ method  to prove the existence of peak functions in
$A(D)$. In fact, they also show how a strengthened form of the
lemma gives greater regularity:

\proclaim 5.4 Theorem. Let $D$ be a smooth bounded pseudoconvex
domain in $\CC$ of finite type. Then there exists a constant
$\alpha
>0$ so that every boundary point of $D$ is a peak point for $A^\alpha(D)$.

\vspace{.1in}

Note that the regularity exponent $\alpha$ is independent of the
boundary point. This exponent can be expressed in terms of the
maximum type of $D$ and a constant (with a similar role to the
constant $C$ in the lemma) that arises in the proof. Their proof
requires H\"older estimates for the $\overline{\partial}
$-equation to pass from local to global peak functions. As we
explain in Section 6 below, Cho later showed how to avoid this
dependence on estimates for the $\overline{\partial} $-equation.

We remark that the relevance of the Bergman kernel to the
existence of peak functions can be seen in the following way. Let
$D$ be a smooth bounded pseudoconvex domain in $\Cn$ with Bergman
kernel $K$. Assume that the kernel can be extended to a function,
still denoted by $K$, on
$(\overline{D}\times\overline{D})\setminus\{(z,z)\colon z\in\b
D\}$ so that, for fixed $p\in \b D$, the function $h(z)= K(z,p)$
is holomorphic on $D$ and continuous  on
$\overline{D}\setminus\{p\}$. (This extension is known to hold on
all pseudoconvex domains of finite type in $\Cn$, where $h$ is
known even to be smooth on $\overline{D}\setminus\{p\}$. For the
case of strictly pseudoconvex domains see the seminal work of
Kerzman in \cite{kerz}. See also the work of Bell in \cite{bell},
Boas in \cite{boas}, and Chen in \cite{chen}.) Then the following
holds.

\proclaim 5.5 Proposition. With the preceding notation, if $h(z)$
tends to infinity as $z$ approaches $p$ through points of $D$,
then $p$ is a local peak point for $A(D)$.

\noindent{\bf Proof:} The function $1/h$ can be extended to a
function in $A(D\cap V)$ (for some neighborhood $V$ of $p$) whose
only zero in the closure of $D\cap V$ is at $p$. We will see in
the last section (Theorem 7.2)  that it follows that $p$ is a
local peak point for $A(D)$. \hfill $\Box$

\vspace{.1in}

We remark that the condition on the kernel in the proposition
should not be confused with the statement that $K(z,z)$ tends to
infinity as $z$ approaches the boundary from the interior. As we
noted above, that statement is valid on all smooth bounded
pseudoconvex domains. The condition on the kernel in the
proposition is known to hold at strictly pseudoconvex boundary
points  and on the relatively few pseudoconvex domains of finite
type for which the Bergman kernel can be computed explicitly.

\vspace{.2in}

{\bf 6. Domains of finite type in $\Cn$}

\vspace{.2in}

The constructions of peak functions described in Sections 4 and 5
can be extended to certain domains in $\Cn$. A recurring theme is
the usefulness of bumping outward the domain near the boundary
point in question.

In \cite{bf1}, Bedford and Forn{\ae}ss show that their results in
$\CC$ imply the following. Suppose that $D$ is a smooth bounded
pseudoconvex domain in $\Cn$ of finite type at $p$ and that the
Levi form at $p$ (considered as a Hermitian form on the complex
tangent space) has corank 1 (that is, it has $n-2$ positive
eigenvalues). Then $p$ is a peak point for $A^\alpha(D)$ for some
$\alpha>0$. Of course, the rank condition is automatically
satisfied in $\CC$.

In \cite{fomc}, Forn{\ae}ss and McNeal show that their method
involving the Bergman kernel also applies to the class of
decoupled pseudoconvex domains of finite type in $\Cn$, defined as
follows. Let $D$ be a smooth bounded pseudoconvex domain of finite
type at $p$. We say that $D$ is decoupled near $p$ if there exist
holomorphic coordinates $z=(z_1,\ldots,z_{n-1},w)$ near $p$ in
which $p=0$ and in which $D$ is defined near $p$ by $$\mbox{Re
}w+\sum_{j=1}^{n-1} \phi_j(z_j)<0.$$ Here each $\phi_j$ is a
smooth subharmonic, but not harmonic, function on ${\bf C}$.
Forn{\ae}ss and McNeal show that if $D$ is a decoupled
pseudoconvex domain of finite type in $\Cn$ then there exists a
constant $\alpha
>0$ so that every boundary point of $D$ is a peak point for
$A^\alpha(D)$. Their method of constructing local peak functions
is exactly the same as for domains of finite type in $\CC$---the
required conditions on the Bergman kernel and the
$\overline{\partial} $-Neumann operator are known to be satisfied
for decoupled domains.

The method of Forn{\ae}ss and McNeal was extended by Cho in
\cite{cho1} and \cite{cho2}. (See also Kol\'a\v{r}'s work in \cite
{kolp}.) Cho uses bumping families he had constructed in earlier
work to avoid the requirement of uniform estimates for the
$\overline{\partial} $-equation. (Cho's bumping fixes the boundary
on a neighborhood of the given point.) His method produces a peak
function at those boundary points of domains of finite type in
$\Cn$ near which the domain is convex. Here he uses work of McNeal
estimating the Bergman kernel on convex domains. The peak function
belongs to $A^\alpha$ (for some $\alpha>0$) and extends
holomorphically across the boundary away from arbitrarily small
neighborhoods of the peak point. Cho's method also applies to the
classes of pseudoconvex domains  of finite type we have already
mentioned: domains for which the Levi form has corank at most 1,
and decoupled domains.

As we noted in Section 3, Diederich and Forn{\ae}ss in
\cite{difoc1} gave an explicit formula for a local strong support
function in $A^\omega$ at each boundary point of a convex domain
of finite type. Their result covers the locally convex case
considered by Cho.

Now we turn to results  that depend on bumping outward the domain
at a given boundary point. This is a different kind of bumping
than Cho's: Here a bumping of $D$ near $p\in\partial D$ is a
domain containing $\overline{D}\setminus\{p\}$ near $p$. In
\cite{bf1}, Bedford and Forn{\ae}ss show how to extend the sector
method to the domain $D$ in $\Cn$ defined by $$\mbox{Re
}w+P(z')<0,$$ where $P$ is a plurisubharmonic polynomial on ${\bf
C}^{n-1}$ that is homogeneous of degree $2k>0$ and satisfies the
following condition: The function $P_\epsilon$ on ${\bf C}^{n-1}$
defined by
$$P_\epsilon(z')=P(z')-\epsilon|z_1|^{2k}-\cdots-\epsilon|z_{n-1}|^{2k}$$
is plurisubharmonic for some $\epsilon>0$. (Recall that an upper
semicontinuous function is plurisubharmonic if it is subharmonic
along each complex line. Such a function $u$ is strictly
plurisubharmonic if  for each point there exists  $\epsilon>0$ so
that the function $u(z)-\epsilon ||z||^2$ is plurisubharmonic near
that point.) Thus the domain defined by $$\mbox{Re
}w+P_\epsilon(z')<0$$ represents a bumping of $D$. Now we explain
why the condition on $P$ and the existence of a bumping are
important for the construction of peak functions by Bedford and
Forn{\ae}ss.

As in the construction in $\CC$, the first step in $\Cn$ is to
intersect the cone $C$ defined by $$\mbox{Re }(w^{2k})+P(z')<0$$
with complex lines through $0$ and to study the resulting sectors.
By assumption, $P$ is strictly subharmonic away from the origin on
each such complex line. This implies that the sectors vary
smoothly with the complex line, so they may be considered as
points of a Riemann domain $S$ over complex projective space of
dimension $n-1$. To obtain the desired function $g$ linear on each
sector, one needs that $S$ is Stein. In \cite{bf1}, Bedford and
Forn{\ae}ss show that the condition on $P$ implies that $S$ is
locally pseudoconvex. Then $S$ is Stein by Fujita's solution of
the Levi problem for Riemann domains over projective space.

So far the existence of a bumping domain has not been used
directly. But the construction of peak functions---in particular
the holomorphic extension across the boundary away from the
origin---requires that $g$ be defined on a conical neighborhood of
$C$. In $\CC$, this is accomplished in \cite{bf1} by examining the
projection of  the boundary of $S$ to complex projective space. In
$\Cn$, under the assumptions on $P$ one constructs a larger
Riemann domain $S_\epsilon$ by applying the sector construction to
the cone $C_\epsilon$ defined by $$\mbox{Re
}(w^{2k})+P_\epsilon(z')<0.$$ Note that $S$ is a relatively
compact subset of $S_\epsilon$ and $C_\epsilon$ is a conical
neighborhood of $C$. The desired function $g$ on $C_\epsilon$ is
then obtained from a section of a line bundle over  $S_\epsilon$.

We remark that, because of the bumping procedure, the construction
of Bedford and Forn{\ae}ss is not restricted to  pseudoconvex
domains with a defining function of the form $\mbox{Re }w+P(z')$.
To explain this, we recall from the beginning of Section 3 the
function $A$ appearing in a standard form for a defining function.
The method of Bedford and Forn{\ae}ss produces a peak function at
every boundary point for which the initial homogeneous term $P$ of
the expansion of $A$ satisfies the condition in their result
(namely, that the associated function $P_\epsilon$ be
plurisubharmonic for some $\epsilon>0$). Thus the defining
function for the general domain has the form $$\mbox{Re
}w+P(z')+O(|\mbox{Im }w|^2+||z'|| |\mbox{Im }w|+||z'||^{2k+1}).$$
We call the domain defined by $\mbox{Re }w+P(z')<0$ a homogeneous
model for such a general domain. The point is that the
higher-order terms in the defining function are controlled by
means of the conical neighborhood.

In \cite{no:pk}, Noell showed how to weaken the assumption on the
leading homogeneous term $P$. He showed that if $P$ is not
harmonic on any complex line through the origin then $0$ is a peak
point for $A^\alpha$ (for some $\alpha>0$) relative to the model
domain, and is in fact a peak point for any smooth bounded
pseudoconvex domain with such a model. The key proposition in
\cite{no:pk} is the following bumping result.

\proclaim 6.1 Proposition. Let $P$ be a homogeneous
plurisubharmonic polynomial on ${\bf C}^{n-1}$. Assume that $P$ is
not harmonic on any complex line through the origin. Then there
exists a function $F$ on ${\bf C}^{n-1}\setminus \{0'\}$ that is
smooth, positive, and homogeneous of the same degree as $P$, and
so that $P-F$ is strictly plurisubharmonic.

\vspace{.1in}

Note that $F$ need not be smooth at the origin.

\vspace{.1in}

\noindent{\bf Sketch of proof:} The proof uses the original
geometric setting: Let $D$ be the homogeneous  model defined by
$$\mbox{Re }w+P(z')<0.$$ The function $F$ is described in terms of
a family of pseudoconvex neighborhoods of
$\overline{D}\setminus\{z'=0\}$ that are strictly pseudoconvex
away from $\{z'=0\}$. These neighborhoods are constructed using a
modification of the technique introduced by Diederich and
Forn{\ae}ss in \cite{difonb}. Their technique requires the
existence of a certain stratification of the boundary, and in the
case of $D$ such a stratification exists by the work of  Diederich
and Forn{\ae}ss in \cite{diforb}.

The strictly pseudoconvex neighborhoods have defining functions of
the form $$\mbox{Re }w+P(z')-\psi(z')||z'||^{2k-2},$$ where $2k$
is the degree of $P$, $\psi$ is homogeneous of degree 2, and away
from the origin $\psi$ is smooth and positive. Once such a
function $\psi$ is constructed to give a domain that is strictly
pseudoconvex away from $\{z'=0\}$, one can define
$F(z')=\psi(z')||z'||^{2k-2}$. To construct $\psi$, note that the
homogeneity of $P$ implies that the projection of each stratum to
${\bf C}^{n-1}$ can be taken to be invariant under dilations.
Locally $\psi$ is defined to be  a large negative multiple of the
distance squared to such a projected stratum, plus $||z'||^2$.
\hfill $\Box$

\vspace{.1in}

The condition on $P$ in Noell's result is equivalent to requiring
that the homogeneous model $$\mbox{Re }w+P(z')<0$$ be pseudoconvex
and of finite type at the origin. His result does not apply to
general pseudoconvex domains of finite type for the following
reason: Even when the function $A$ (from the standard form for a
defining function) does not vanish on any nontrivial complex
variety, the leading homogeneous term $P$ of $A$ may vanish on a
complex line through the origin. For example, define $A$  on $\CC$
by $A(z_1,z_2)=|z_1|^2+|z_2|^4$. Then $P(z_1,z_2)=|z_1|^2$, which
vanishes on the complex line $z_1=0$. In this example  assigning
weight 2 to $z_1$ and weight 1 to $z_2$ will avoid this
difficulty. In \cite {herb}, Herbort showed how to construct a
peak function on certain weighted homogeneous model domains. His
method uses H{\"o}rmander's theory of solving the
$\overline{\partial} $-equation in $L^2$ spaces with weights, as
in the construction of peak functions by Forn{\ae}ss and Sibony.

To discuss the general weighted homogeneous case we use the
terminology of Yu  in \cite{yupk} (see also \cite{yudi},
\cite{yuco}, and \cite{yukr}), who proved the definitive results
along these lines. Quite similar results were proved by Diederich
and Herbort in \cite{dihe}. A multiindex
$\Lambda=(\lambda_1,\ldots,\lambda_m)$ is called a multiweight if
$1\geq \lambda_1\geq\cdots\geq\lambda_m>0$. A function $h$ on
${\bf C}^m$ is said to be $\Lambda$-homogeneous if
$$h(t^{\lambda_1}z_1,\ldots,
t^{\lambda_m}z_m)=th(z_1,\ldots,z_m)$$ for every $t>0$ and every
$(z_1,\ldots,z_m)\in {\bf C}^m$. Let $P$ be a  real-valued
function on ${\bf C}^{n-1}$ that is  smooth on ${\bf
C}^{n-1}\setminus \{0'\}$ , and define $$D=\{(z',w)\in\Cn\colon
\mbox{Re }w+P(z')<0\}.$$ We say that $D$ is a
$\Lambda$-homogeneous model if $P$ is $\Lambda$-homogeneous and
plurisubharmonic but not pluriharmonic (i.e., not the real part of
a holomorphic function), and we say that $D$ is a
$\Lambda$-homogeneous polynomial model if in addition $P$ is a
polynomial. If $D$ is a $\Lambda$-homogeneous model with function
$P$, then $D$ is called $h$-extendible if there exists a
$\Lambda$-homogeneous function $F$ on ${\bf C}^{n-1}$ that on
${\bf C}^{n-1}\setminus \{0'\}$ is of class $C^1$ and positive,
and so that $P-F$ is plurisubharmonic. (Compare this definition
with the conclusion of Proposition 6.1.) The following bumping
result is proved in \cite{yupk} (see also \cite{dihe}).

\proclaim 6.2 Theorem. Every $\Lambda$-homogeneous polynomial
model of finite type is $h$-extendible.

The bumping function $F$ can be taken to be smooth away from the
origin. The proof is formally similar to that of Proposition 6.1,
but it is much more difficult and depends on results of Catlin in
\cite{catl}.

\vspace{.1in}

Earlier, in \cite{yudi} (see also \cite{yupk}), Yu had showed how
to construct peak functions on $h$-extendible models:

\proclaim 6.3 Theorem. If $D$ is an  $h$-extendible model, then
$0$ is a peak point for $A(D)$.

The method uses H{\"o}rmander's theory of solving the
$\overline{\partial} $-equation in $L^2$ spaces with weights, as
in the construction of peak functions by Forn{\ae}ss and Sibony.

By combining Theorem 6.2 and Theorem 6.3, one gets a peak function
at the origin (extending holomorphically across the boundary away
from the origin) for every $\Lambda$-homogeneous polynomial model
of finite type.

 To explain how these results apply to domains more
general than models, we need more definitions. Fix a boundary
point  $p$ of a smooth bounded pseudoconvex domain in $\Cn$. For
$1\leq q\leq n$ we denote by $\Delta_q(p)$ the maximum order of
contact of $q$-dimensional complex analytic varieties with the
boundary at $p$, the D'Angelo $q$-type. Now assume that $p$ is a
point of finite type (so $\Delta_1(p)$ is finite). In \cite{catl},
Catlin defined the multitype at $p$. The actual definition is
somewhat complicated, but briefly the multitype is an $n$-tuple
$(m_1(p), \ldots , m_n(p))$ of rational numbers determined by
measuring orders of vanishing of a defining function in the
coordinate direction $z_j$ with the weight $m_j(p)$ attached.
Always $m_1(p)=1$, $m_2(p)=\Delta_{n-1}(p)$, and $m_1(p)\leq
m_2(p)\leq \ldots \leq m_n(p)$. Further, if the Levi form has rank
$q$ at $p$ then $m_j(p)=2$ when $2\leq j\leq q+1$ and $m_j(p)>2$
when $j>q+1$. Catlin proved that if $1\leq q\leq n$ then
$$m_{n+1-q}(p)\leq\Delta_{q}(p).$$ Here is the main result of
\cite{dihe} and \cite{yupk} with regard to peak functions on
bounded domains.

\proclaim 6.4 Theorem. Let $D$ be a smooth bounded pseudoconvex
domain in $\Cn$ and $p$ a boundary point of finite type. If
$$m_{n+1-q}(p)=\Delta_{q}(p)$$ when $1\leq q\leq n$, then $p$ is a
peak point for $A(D)$.

In \cite{dihe}, in addition it is proved that there exists a local
peak function in $A^\alpha$ for some $\alpha>0$, so there is a
global peak function in some $A^\alpha$ if, say, $D$ is of finite
type.

The proof of the theorem depends on relating the domain to a model
domain. To explain this we use the notation of the theorem and
define the multiweight $\Lambda=(1/m_2(p),\ldots,1/m_n(p))$. If
$r$ is a defining function near $p$, then Catlin's work gives
local holomorphic coordinates $(z',w)$ in which $p=0$ and
$$r(z',w)=\mbox{Re }w+P(z')+R(z',w),$$ where $P$ is a
$\Lambda$-homogeneous plurisubharmonic polynomial with no
pluriharmonic terms and $R$ is a smooth function of strictly
higher order than $P$ (in terms of the coordinates weighted
according to Catlin's multitype). Then $$\mbox{Re }w+P(z')<0$$
defines a $\Lambda$-homogeneous polynomial model associated to $D$
at $p$.  We say that $D$ is $h$-extendible at $p$ if it has an
associated $\Lambda$-homogeneous polynomial model that is
$h$-extendible. Yu proved that in this case there exists an
$h$-extendible model $D'$ so that $\overline D\setminus\{p\}$ is
contained in $D'$ near $p$. Now we can state the  remarkable
characterization of $h$-extendibility given in \cite{dihe} and
\cite{yupk}.

\proclaim 6.5 Theorem. Let $D$ be a smooth bounded pseudoconvex
domain in $\Cn$ and $p$ a boundary point of finite type. Then $D$
is $h$-extendible at $p$ if and only if
$$m_{n+1-q}(p)=\Delta_{q}(p)$$ when $1\leq q\leq n$.

\vspace{.1in}

We omit the proof. Now the proof of Theorem 6.4 is clear: Apply
Theorem 6.5 to see that $D$ is $h$-extendible at $p$, and use Yu's
result to get an $h$-extendible model $D'$ so that $\overline
D\setminus\{p\}$ is contained in $D'$ near $p$. Then apply the
earlier construction of peak functions on $\Lambda$-homogeneous
polynomial models of finite type to get a peak function relative
to $D'$ at $p$. This function is also a local peak function
relative to $D$ and extends holomorphically across the boundary
away from $p$.

We note that Theorem 6.4 implies all of the results presented in
this section on the peak function problem: All of the classes of
bounded domains considered in this section---domains for which the
Levi form has corank at most 1, decoupled domains, and convex
domains---are known to be $h$-extendible.

Recently, Bharali and Stens{\o}nes have proved bumping results for
certain domains in ${\bf C}^3$ that are not $h$-extendible. See
\cite{bhst}.

\vspace{.2in}

{\bf 7. Miscellany}

\vspace{.2in}

In this section we list some miscellaneous results related to the
peak point problem and its generalizations.

One natural question is whether peak functions can be chosen to
vary in a regular way from point to point. In \cite{grah}, Graham
showed that on smooth bounded strictly pseudoconvex domains one
can find peak functions in $A^\omega$ that vary continuously with
the boundary point.
    Forn{\ae}ss and Krantz proved a very general result in
\cite{fokr}: If $X$ is a compact metric space and ${\cal A}$ is a
closed subalgebra of $C(X)$ whose set of peak points is ${\cal
P}$, then there is a continuous function $\Phi\colon{\cal
P}\rightarrow{\cal A}$ such that $\Phi(x)$ is a peak function for
$x$ when $x\in{\cal P}$.

One important result on peak points proved by Gamelin in
\cite{game} involves an approach very different from any
considered so far. Recall that a subset $K$ of $\Cn$ is said to be
circled if $(\lambda_1 z_1, \ldots, \lambda_n z_n)\in K$ whenever
$(z_1, \ldots, z_n)\in K$ and $\lambda_1 , \ldots, \lambda_n$ are
complex numbers of modulus 1. Let $K$ be a compact, connected,
circled subset of $\Cn$. We denote by $H(K)$ the uniform closure
in $C(K)$ of the functions holomorphic in a neighborhood of $K$.
The rational convex hull of $K$, denoted $r(K)$, is the set of
points $z\in\Cn$ with the property that every holomorphic
polynomial that vanishes at $z$ also has a zero on $K$. (For such
$K$ the set $r(K)$ may be identified with the spectrum of $H(K)$.)
Here is the main result from \cite{game}.

\proclaim 7.1 Theorem.  Let $K$ be a compact, connected, circled
subset of $\Cn$. If $p\in K$ is not a peak point for $H(K)$, then
$p$ belongs to an analytic disc that is contained in $r(K)$.

The proof uses methods from the theory of function algebras.

\vspace{0.1in}

In the introduction we noted that analytic structure in the
boundary is an obstruction to the existence of peak functions. In
\cite{yuce}, Yu showed that even without such analytic structure a
peak function may fail to exist: He constructed a smooth bounded
pseudoconvex domain $D$ in ${\bf C}^3$ so that $\partial D$
contains no nontrivial analytic set and so that some boundary
point is not a local peak point for $A(D)$. The domain constructed
by Yu is not of finite type, but it is B-regular, which means
roughly that there is no pluripotential theory on the boundary.
(See Sibony's survey paper \cite{sib:asp} for a discussion of
B-regularity.) An  example of a related phenomenon was given by
Noell in  \cite{no:inty}: There is a smooth bounded convex domain
$D$ in $\CC$ so that $D$ is strictly pseudoconvex except along a
line segment and so that, for all $\alpha >0$, each point of that
segment is not a peak point for $A^\alpha(D)$. The domain is of
infinite type along the line segment, and every point of the
segment is a peak point for $A(D)$.

The notion of a peak point can be generalized as follows. If $D$
is a domain in $\Cn$ with smooth boundary,  a compact subset $K$
of $\partial D$ is a peak set relative to $D$ for a space ${\cal
A}$ if there exists a function $f\in {\cal A}$ so that $f=1$ on
$K$ and $|f|<1$ on $\overline{D}\setminus K$. As for peak points,
we use the terms peak function and strong support function. Peak
sets have been the subject of extensive study, but even in the
case $n=1$ a complete characterization of peak sets for $A^\alpha$
is unknown if $0<\alpha<1$.

We state here an open question that is a natural generalization of
the peak point problem.

\proclaim Question. If $D$ is a smooth bounded pseudoconvex domain
in $\Cn$, does every boundary point of $D$ belong to a peak set
for $A(D)$?

An affirmative answer would imply completeness in the standard
invariant metrics.  We remark that clearly every boundary point of
a bounded convex domain belongs to a peak set for $A^\omega$.

A result from the theory of function algebras implies that if $D$
is a bounded domain then every peak set for $A(D)$ contains a peak
point for $A(D)$. (See Corollary 2.4.6 of \cite{brow}.) In the
example of Noell mentioned above, the line segment is a peak set
for $A^\omega$. This illustrates the fact that the  result from
function algebra theory does not extend to $A^\alpha$ for
$\alpha>0$.

Sets more general than peak sets have also been studied.  If $D$
is a domain in $\Cn$ with smooth boundary,  a compact subset $K$
of $\partial D$ is a zero set relative to $D$ for a space ${\cal
A}$ if there exists a function $f\in {\cal A}$ so that $f=0$ on
$K$ and $f\neq 0$ on $\overline{D}\setminus K$. We call $f$ a
support function. Of course, every peak set is a zero set, and
each strong support function is a support function.

The original paper \cite{bf1} of Bedford and Forn{\ae}ss
constructed a support function in $A^\infty$ at each boundary
point of pseudoconvex domains of finite type in $\CC$. The support
function they construct vanishes to infinite order at the point in
question. As noted in \cite{bf1}, the argument of  Kohn and
Nirenberg in \cite{koni} shows that for the Kohn-Nirenberg domain
every support function for the origin must vanish to infinite
order there.

At the other extreme from the case of infinite order of vanishing,
we have the case when, for a point $p$ in the boundary of a domain
$D$, there is a local support function $f\in A^\omega(D)$ with
nonzero gradient at $p$. Then it follows that there is a complex
submanifold $X$ of complex codimension 1 in a neighborhood $V$ of
$p$ so that $X\cap\overline{D}\cap V=\{p\}$. We call such a set
$X$ a local support manifold for $D$ at $p$. We remark that, by
 the holomorphic Hopf lemma (2.6), if $p$ is a
local peak point for $A^\omega$ then there is a local support
manifold at $p$. We refer to  Th\'eor\`eme 1 (and its proof) in
the paper \cite{hasib:qc} of Hakim and Sibony for a very useful
generalization of this fact that covers the case when $p$ belongs
to a peak set for $A^\infty(D)$, with $D$ of finite type at $p$.

The existence of a local support manifold is closely related to
the strict type conditions considered in Section 3. In fact,
Range's strict type condition (the weakest of the three conditions
considered) at a point clearly implies the existence of a local
support manifold there. For the converse, in \cite{bloo} Bloom
observed that if there is a local support manifold for $p$, and if
the boundary is real-analytic near $p$, then Range's strict type
condition holds at $p$.

The example of Bloom that we considered in Proposition 3.6 has a
local support manifold at the origin (because the domain is of
strict type in the sense of Range at the origin). This example
shows that the existence of a smooth peak function does not follow
from the existence of a local support manifold. It is true,
however, that if there is a local support manifold at a boundary
point $p$ of a smooth bounded domain $D$ then $p$ is a local peak
point for $A(D)$. (If also $D$ is pseudoconvex then $p$ is a
global peak point for $A(D)$.)  This was proved by Hakim and
Sibony in \cite{hasib:qc} and by Range in \cite{range2}. One
definitive result along these lines was proved by Verdera in
\cite{ver}:

\proclaim 7.2 Theorem. Let $D$ be a smooth bounded pseudoconvex
domain in $\Cn$. If $K\subset\partial D$ is a zero set for $A(D)$
then $K$ is a peak set for $A(D)$.

The elegant proof uses results such as the identification of the
spectrum of $A(D)$ and the Arens-Royden Theorem. In the case of
strictly pseudoconvex domains this theorem was proved by Chaumat
and Chollet in \cite{chch} and by Weinstock in \cite{wei}. Here is
a proof (from \cite{chch}) for the case when $D$ is simply
connected: Let $D$ be a simply connected domain and $K\subset
\partial D$ a zero set for $A(D)$, with support function $f$. We
may assume that $|f|<1$ on $\overline{D}$. Then there exists a
function $\phi$ holomorphic on $D$ so that $\exp{\phi}=f$, and we
have $\mbox{Re }\phi<0$ on $D$. It is easy to see that $-1/\phi$
extends to a strong support function for $K$ in $A(D)$. Hence $K$
is a peak set for $A(D)$. We remark that as a special case we have
that if $D$ is any smooth bounded domain with $p\in\partial D$,
and if $\{p\}$ is a local zero set for $A(D)$, then $p$ is a local
peak point for $A(D)$. This local result is all we needed in the
proof of Proposition 5.5.

We also mention the notion of a plurisubharmonic peak function for
a domain $D$ at a point $p\in\partial D$, namely, a function
$\psi\in C(\overline{D})$ that is plurisubharmonic on $D$ and
satisfies $\psi(p)=0$ and $\psi<0$ on
$\overline{D}\setminus\{p\}$. Sibony showed (see Theorem 2.3 of
\cite{sib:asp}) that a smooth bounded pseudoconvex domain in $\Cn$
is B-regular if and only if there is a plurisubharmonic peak
function at each boundary point. Clearly every peak point for
$A(D)$ has a plurisubharmonic peak function, but in light of Yu's
example above the converse is false. The construction of peak
functions by Forn{\ae}ss and Sibony in \cite{fosib} does suggest a
connection between plurisubharmonic peak functions and
(holomorphic) peak functions in certain settings: One of the main
results of \cite{fosib} is the existence, on each pseudoconvex
domain of finite type in $\CC$, of a family of plurisubharmonic
peak functions (one for each  for boundary point) satisfying
optimal estimates and with uniform control on regularity. (They
prove a similar result on convex domains of finite type in $\Cn$.)
We remark that for some applications the existence of
plurisubharmonic peak functions is sufficient. In \cite{fosib},
for example, Forn{\ae}ss and Sibony show how to deduce optimal
subelliptic estimates for the $\overline{\partial} $-Neumann
problem from the family described above. See also Herbort's work
in \cite{herb2} connecting plurisubharmonic peak functions to
estimates for the Bergman kernel and metric.

\vspace{.2in}

{\sc Department of Mathematics, Oklahoma State University,
Stillwater, OK} \ 74078

{\em E-mail address:} {\tt noell@math.okstate.edu}


\vspace{.2in}

\begin{thebibliography}{9}

\bibitem{ba} R. Basener, {\em Peak points, barriers and pseudoconvex boundary
points\/},  Proc.\ Amer.\ Math.\ Soc.\  {\bf 65}  (1977), no.\ 1,
89--92.

\bibitem{bf1} E. Bedford and J. E. Forn{\ae}ss, {\em A construction of peak functions
on weakly pseudoconvex domains\/},  Ann.\ of Math.\ (2)  {\bf 107}
(1978), no.\ 3, 555--568.

\bibitem{bf2} E. Bedford and J. E. Forn{\ae}ss, {\em
Biholomorphic maps of weakly pseudoconvex domains\/},  Duke Math.\
J. {\bf 45}  (1978), no.\ 4, 711--719.

\bibitem{bell} S. Bell, {\em Differentiability of the Bergman kernel and pseudolocal
estimates\/}, Math.\ Z. {\bf 192} (1986), no.\ 3, 467--472.

\bibitem{bhst} G. Bharali and B. Stens{\o}nes, {\em Plurisubharmonic polynomials and bumping\/}, arXiv:0709.3993v2,
to appear in Math.\ Z.


\bibitem{bloo} T. Bloom, {\em ${\cal C}\sp{\infty }$ peak functions
for pseudoconvex domains of strict type\/},  Duke Math.\ J. {\bf
45} (1978), no.\ 1, 133--147.

\bibitem{boas} H. Boas, {\em Extension of Kerzman's theorem on differentiability of the Bergman
kernel function\/}, Indiana Univ.\ Math.\ J. {\bf 36} (1987), no.\
3, 495--499.

\bibitem{brow} A. Browder, {\em Introduction to function algebras\/}, W. A. Benjamin,
Inc., New York-Amsterdam, 1969.

\bibitem{catl} D. Catlin, {\em Boundary invariants of pseudoconvex
domains\/}, Ann.\ of Math.\ (2) {\bf 120} (1984), no.\ 3,
529--586.

\bibitem{chch} J. Chaumat and A.-M. Chollet, {\em Ensemble de z\'eros,
ensembles pics et d'interpolation pour $A(D)$\/}, Colloque
d'Analyse Harmonique et Complexe, 5 pp.\ (not consecutively
paged), Univ.\ Aix-Marseille I, Marseille, 1977.

\bibitem{chen} S. C. Chen, {\em A counterexample to the
differentiability of the Bergman kernel function\/}, Proc.\ Amer.\
Math.\ Soc.\ {\bf 124} (1996), no.\ 6, 1807--1810.

\bibitem{cho1} S. Cho, {\em A construction of peak functions on locally convex domains in $
{\bf C}\sp n$\/},  Nagoya Math.\ J. {\bf  140}  (1995), 167--176.

\bibitem{cho2} S. Cho, {\em Peak function and its applications\/},  J. Korean
Math.\ Soc.\ {\bf 33}  (1996),  no.\ 2, 399--411.

\bibitem{dabo} J. D'Angelo, {\em Several complex variables and the geometry of real
hypersurfaces\/}, Studies in Advanced Mathematics, CRC Press, Boca
Raton, FL, 1993.

\bibitem{difo} K. Diederich and J. E. Forn{\ae}ss, {\em Pseudoconvex domains: an
example with nontrivial Nebenh\"ulle\/},  Math.\ Ann.\  {\bf 225}
(1977), no.\ 3, 275--292.

\bibitem{difonb} K. Diederich and J. E. Forn{\ae}ss, {\em Pseudoconvex
domains: existence of Stein neighborhoods\/},  Duke Math.\ J. {\bf
44} (1977), no.\ 3, 641--662.

\bibitem{diforb} K. Diederich and J. E. Forn{\ae}ss, {\em
 Pseudoconvex domains with real-analytic boundary\/},  Ann.\ Math.\ (2)
{\bf 107}  (1978), no.\ 2, 371--384.

\bibitem{difoc1} K. Diederich and J. E. Forn{\ae}ss,
{\em Support functions for convex domains of finite type\/},
Math.\ Z. {\bf  230}  (1999),  no.\ 1, 145--164.

\bibitem{difoc2} K. Diederich and J. E. Forn{\ae}ss,
{\em  Lineally convex domains of finite type: holomorphic support
functions\/},  Manuscripta Math.\ {\bf  112}  (2003),  no.\ 4,
403--431.

\bibitem{dihe} K. Diederich and G. Herbort, {\em Pseudoconvex
domains of semiregular type\/}, Contributions to complex analysis
and analytic geometry, 127--161, Aspects Math., E26, Vieweg,
Braunschweig, 1994.

\bibitem{fo:c1} J. E. Forn{\ae}ss, {\em Peak points on weakly pseudoconvex domains\/},
Math.\ Ann.\  {\bf 227}  (1977), no.\ 2, 173--175.

\bibitem{fo:dbar} J. E. Forn{\ae}ss, {\em Sup-norm estimates for $\overline{\partial}$
in $\CC$\/},  Ann. of Math. (2)  {\bf 123}  (1986),  no.\ 2,
335--345.

\bibitem{fokr} J. E. Forn{\ae}ss and S. G. Krantz, {\em Continuously
varying peaking functions\/},  Pacific J. Math.\ {\bf  83} (1979),
no.\ 2, 341--347.

\bibitem{fomc} J. E. Forn{\ae}ss and J. McNeal, {\em A construction of peak
functions on some finite type domains\/}, Amer.\ J. Math.\ {\bf
116} (1994),  no.\ 3, 737--755.

\bibitem{fore} J. E. Forn{\ae}ss and C. Rea, {\em Local holomorphic extendability
and nonextendability of CR-functions on smooth boundaries\/},
Ann.\ Scuola Norm.\ Sup.\ Pisa Cl.\ Sci.\ (4)  {\bf 12}  (1985),
no.\ 3, 491--502.

\bibitem{fosib} J. E. Forn{\ae}ss and N. Sibony, {\em Construction of P.S.H. functions
 on weakly pseudoconvex domains\/},  Duke Math.\ J. {\bf 58}  (1989),  no.\ 3,
 633--655.

\bibitem{fost} J. E. Forn{\ae}ss and B. Stens{\o}nes,
{\em Lectures on counterexamples in several complex variables\/},
Math.\ Notes, 33, Princeton Univ.\ Press, Princeton, NJ, 1987.

\bibitem{free} M. Freeman, {\em Local complex foliation of real
submanifolds\/}, Math.\ Ann.\  {\bf 209}  (1974), 1--30.

\bibitem{game} T. W. Gamelin, {\em Peak points for algebras on circled sets\/},
Math.\ Ann.\  {\bf 238}  (1978), no.\ 2, 131--139.

\bibitem{gase} R. Gay and A. Sebbar, {\em  Division et extension dans l'alg\`ebre
$A\sp \infty(\Omega)$ d'un ouvert pseudo-convexe \`a bord lisse de
$\Cn$\/},  Math.\ Z. {\bf 189}  (1985),  no.\ 3, 421--447.

\bibitem{grah} I. Graham, {\em Boundary behavior of the Carath\'eodory and Kobayashi
metrics on strongly pseudoconvex domains in ${\bf C}\sp{n}$ with
smooth boundary\/},  Trans.\ Amer.\ Math.\ Soc.\ {\bf 207}
(1975), 219--240.

\bibitem{hasib:fs} M. Hakim and N. Sibony, {\em Fronti\`ere de \^Silov et spectre de $A(\overline{D})$
pour des domaines faiblement pseudoconvexes\/}, C. R. Acad.\ Sci.\
Paris S\'er.\ A-B  {\bf 281}  (1975), no.\ 22, Aii, A959--A962.

\bibitem{hasib:qc} M. Hakim and N. Sibony, {\em Quelques conditions pour l'existence de fonctions pics dans des
domaines pseudoconvexes\/},  Duke Math.\ J.  {\bf 44}  (1977),
no.\ 2, 399--406.

\bibitem{herb} G. Herbort, {\em Invariant metrics and peak functions on
pseudoconvex domains of homogeneous finite diagonal type\/},
Math.\ Z. {\bf 209} (1992), no.\ 2, 223--243.

\bibitem{herb2} G. Herbort, {\em Localization lemmas for the Bergman metric at
plurisubharmonic peak points\/},  Nagoya Math.\ J. {\bf 171}
(2003), 107--125.

\bibitem{iord} A. Iordan, {\em Pseudoconvex domains with peak functions at each point
of the boundary\/},  Pacific J. Math.\ {\bf  133}  (1988),  no.\
2, 277--287.

\bibitem{kerz} N. Kerzman, {\em  The Bergman kernel function. Differentiability
at the boundary\/},  Math.\ Ann.\ {\bf 195} (1972), 149--158.

\bibitem{ko:st} J. J. Kohn, {\em Boundary behavior of $\overline{\partial} $ on weakly pseudo-convex manifolds
of dimension two\/}, J. Differential Geometry {\bf 6} (1972),
523--542.

\bibitem{ko:gl} J. J. Kohn, {\em Global regularity for $\overline{\partial} $ on weakly pseudo-convex manifolds\/},
 Trans.\ Amer.\ Math.\ Soc.\  {\bf 181}  (1973), 273--292.

\bibitem{ko:gl2} J. J. Kohn, {\em Methods of partial differential equations in complex analysis\/},
Several complex variables (Proc.\ Sympos.\ Pure Math., Vol.\ XXX,
Part 1, Williams Coll., Williamstown, Mass., 1975),  pp. 215--237.
Amer.\ Math.\ Soc., Providence, R.I., 1977.

\bibitem{koni} J. J. Kohn and L. Nirenberg, {\em A pseudo-convex domain not admitting a
holomorphic support function\/},  Math.\ Ann.\  {\bf 201}  (1973),
265--268.

\bibitem{kol} M. Kol\'a\v{r}, {\em Convexifiability and supporting functions in
$\CC$\/},  Math.\ Res.\ Lett.\  {\bf 2}  (1995),  no.\ 4,
505--513.

\bibitem{kolp} M. Kol\'a\v{r}, {\em Peak functions on convex domains\/}, The Proceedings of the 19th
Winter School ``Geometry and Physics'' (Srn\'i, 1999).  Rend.\
Circ.\ Mat.\ Palermo (2) Suppl.\  No.\ 63  (2000), 103--112.

\bibitem{la} G. Laszlo, {\em Peak functions on finite type domains
in $\CC$\/}. Dissertation, Oklahoma State University, 2000.

\bibitem{lev1} E. E. Levi, {\em Studii sui punti singolari essenziali delle funzioni analitiche
di due o pi\`u variabili complesse\/},  Ann.\ Mat.\ Pura Appl.\
{\bf 17} (1910), 61-–87.

\bibitem{lev2} E. E. Levi, {\em Sulle ipersuperficie dello spazio a 4
dimensioni che possono essere frontiera del campo di esistenza di
una funzione analitica di due variabili complesse\/},  Ann.\ Mat.\
Pura Appl.\ {\bf  18}  (1911),   69–-79.

\bibitem{mc} J. McNeal, {\em Lower bounds on the Bergman metric near a point of finite
type\/},  Ann.\ Math.\ (2) {\bf 136}  (1992), no.\ 2, 339--360.

\bibitem{no:peak} A. Noell, {\em Properties of peak sets in
weakly pseudoconvex boundaries in ${\bf C}^2$\/}, Math.\ Z. {\bf
186} (1984), 99--116.

\bibitem{no:inty} A. Noell, {\em Peak points in boundaries not of finite type\/},  Pacific J.
Math.\ {\bf 123}  (1986),  no.\ 2, 385--390.

\bibitem{no:ct} A. Noell, {\em  Interpolation from curves in pseudoconvex
boundaries\/}, Michigan Math.\ J. {\bf 37} (1990), no.\ 2,
275--281.

\bibitem{no:pk} A. Noell, {\em Peak functions for pseudoconvex domains in ${\bf
C}^n$\/}, Several complex variables (Stockholm, 1987/1988),
529--541, Math.\ Notes, 38, Princeton Univ.\ Press, Princeton, NJ,
1993.

\bibitem{no:st} A. Noell and B. Stens{\o}nes, {\em Proper holomorphic maps from weakly pseudoconvex
domains\/},  Duke Math.\ J. {\bf 60}  (1990),  no.\ 2, 363--388.

\bibitem{pf} P. Pflug, {\em \"Uber polynomiale Funktionen auf
Holomorphiegebieten},  Math.\ Z.  {\bf 139}  (1974), 133--139.

\bibitem{pf2} P. Pflug, {\em Quadratintegrable holomorphe Funktionen und die Serre-Vermutung},
Math.\ Ann.\  {\bf 216}  (1975), 285--288.


\bibitem{range1} R. M. Range,  {\em  H\"older estimates for $\overline \partial $ on
convex domains in ${\bf C}\sp{2}$ with real analytic boundary\/},
Several complex variables (Proc.\ Sympos.\ Pure Math., Vol. XXX,
Part 2, Williams Coll., Williamstown, Mass., 1975),  pp. 31--33.
Amer.\ Math.\ Soc., Providence, R.I., 1977.

\bibitem{range2} R. M. Range,  {\em The Carath\'eodory metric and holomorphic maps on a class of weakly
pseudoconvex domains\/},  Pacific J. Math.\ {\bf  78}  (1978),
no.\ 1, 173--189.

\bibitem{range3} R. M. Range,  {\em On H\"older estimates for $\overline \partial u=f$ on weakly pseudoconvex
domains\/}, Several complex variables (Cortona, 1976/1977),  pp.\
247--267, Scuola Norm.\ Sup.\ Pisa, Pisa, 1978.

\bibitem{ross} H. Rossi, {\em Holomorphically convex sets in several complex
variables\/},  Ann.\ of Math.\ (2)  {\bf 74}  (1961), 470--493.

\bibitem{sib:dbar} N. Sibony, {\em Un exemple de domaine pseudoconvexe r\'egulier o\`u l'\'equation $\bar
\partial u=f$ n'admet pas de solution born\'ee pour $f$
born\'ee\/}, Invent.\ Math.\ {\bf 62} (1980/81), no.\ 2, 235--242.

\bibitem{sib:asp} N. Sibony, {\em Some aspects of weakly pseudoconvex domains\/},  Several complex
variables and complex geometry, Part 1 (Santa Cruz, CA, 1989),
199--231, Proc.\ Sympos.\ Pure Math., 52, Part 1, Amer.\ Math.\
Soc., Providence, RI, 1991.

\bibitem{ver} J. Verdera, {\em A remark on zero and peak sets on weakly
pseudoconvex domains\/}, Bull.\ London Math.\ Soc.\ {\bf 16}
(1984), no.\ 4, 411--412.

\bibitem{wei} B. Weinstock, {\em Zero-sets of continuous holomorphic
functions on the boundary of a strongly pseudoconvex domain\/}, J.
London Math.\ Soc.\ (2) {\bf 18} (1978), no.\ 3, 484--488.

\bibitem{yudi} J. Yu, {\em Geometric analysis on weakly pseudoconvex domains\/}. Dissertation, Washington
University, St.\ Louis, 1993.

\bibitem{yuco} J. Yu, {\em Multitypes of convex domains\/},  Indiana Univ.\ Math.\ J. {\bf 41}  (1992),
no.\ 3, 837--849.

\bibitem{yupk} J. Yu, {\em Peak functions on weakly pseudoconvex
domains\/},  Indiana Univ.\ Math.\ J. {\bf 43} (1994), no.\ 4,
1271--1295.

\bibitem{yukr} J. Yu, {\em Weighted boundary limits of the generalized Kobayashi-Royden
metrics on weakly pseudoconvex domains\/},  Trans.\ Amer.\ Math.\
Soc.\ {\bf 347}  (1995),  no.\ 2, 587--614.

\bibitem{yuce} J. Yu, {\em A counterexample to the existence of peaking functions\/},
Proc.\ Amer.\ Math.\ Soc.\ {\bf 125}  (1997),  no.\ 8, 2385--2390.

\end{thebibliography}
\end{document}